\documentclass[12pt]{amsart}
\usepackage{amscd,amssymb}
 \newlength{\baseunit}               
\newcommand{\bpf}{\noindent {\em Proof. }}
\newcommand{\epf}{\qed \vspace{+10pt}}

\newtheorem{tm}{Theorem}[section]
\newtheorem{pr}[tm]{Proposition}
\newtheorem{lm}[tm]{Lemma}

\newtheorem{conj}[tm]{Conjecture}
\newtheorem{defn}[tm]{Definition}

\newcommand{\proj}{\mathbb P}

\newcommand{\eff}{\mathbb F}
\newcommand{\fn}{{\mathbb F}_n}

\newcommand{\oh}{{\mathcal{O}}}

\newcommand{\tth}{\text{th}}
\newcommand{\ti}{\text{irr}}

\newcommand{\al}{\alpha}
\newcommand{\be}{\beta}
\newcommand{\ga}{\gamma}
\newcommand{\Ga}{\Gamma}
\newcommand{\de}{\delta}

\newcommand{\si}{\sigma}
\newcommand{\Up}{\Upsilon}
\newcommand{\wdg}{W^{d,g}}
\newcommand{\vdg}{V^{d,g}}
\newcommand{\ab}{(\al,\be)} 
\newcommand{\abs}{(\al,\be,s)} 
\newcommand{\abG}{(\al,\be,\Ga)} 

\newcommand{\lcm}{\operatorname{lcm}}

\newcommand{\idim}{\operatorname{idim}}

\newcommand{\mbar}{\overline{M}}

\begin{document}
\title{Genus $g$ Gromov-Witten invariants of Del Pezzo surfaces:  
Counting plane curves with fixed multiple points} 
\author{Ravi Vakil}
\date{\today} 
\begin{abstract}
As another application of the degeneration methods of [V3], we count the
number of irreducible degree $d$ geometric genus $g$ plane curves, with
fixed multiple points on a conic $E$, not containing $E$, through an
appropriate number of general points in the plane.  As a special case, we
count the number of irreducible genus $g$ curves in any divisor class $D$
on the blow-up of the plane at up to five points (no three collinear).  We
then show that these numbers give the genus $g$ Gromov-Witten invariants of
the surface.  Finally, we suggest a direction from which the remaining del
Pezzo surfaces can be approached, and give a conjectural algorithm to
compute the genus g Gromov-Witten invariants of the cubic surface.
\end{abstract}
\maketitle
\pagestyle{plain}
\tableofcontents

\section{Introduction}

In this note, we count the number of irreducible degree $d$ geometric
genus $g$ plane curves, with fixed multiple points on a conic $E$, not
containing $E$, through an appropriate number of general points in the
plane.  As a special case, we count the number of irreducible genus $g$
curves in any divisor class $D$ on the blow-up of the plane at up to five
points (no three collinear).  These numbers are the genus $g$ Gromov-Witten
invariants of the surface (Subsection \ref{gwenumerative}).  The genus $g$
Gromov-Witten invariants of $\proj^2$ were already computed in [R] and
[CH3], and those of $\proj^1 \times \proj^1$ and the blow-up of $\proj^2$
at a point ($\eff_1$) were computed in [V3].

Such classical enumerative questions have recently been the object of study
by many people.  Ideas from mathematical physics (cf. the inspiring [KM]
and [DI]) have yielded formulas when $g=0$ on $\proj^2$ (via associativity
relations in quantum cohomology).  Z. Ran solved the analogous
(enumerative) problem for curves of arbitrary genus on $\proj^2$ by
degenerations methods (cf. [R]), and Caporaso and Harris gave a second
solution by different degeneration methods (cf. [CH3]).  These numbers for
irreducible curves are also the genus $g$ Gromov-Witten invariants of
$\proj^2$ (Subsection \ref{gwenumerative}).  P. Di Francesco and
C. Itzykson calculated the genus 0 Gromov-Witten invariants of the plane
blown up at up to six points in [DI], Subsection 3.3.  Y. Ruan and G. Tian
gave recursive formulas for the genus 0 Gromov-Witten invariants of Fano
surfaces, and indicated their enumerative significance ([RT] Section 10).
L. G\"{o}ttsche and R. Pandharipande later derived recursive formulas
for the genus 0 Gromov-Witten invariants of the plane blown up at any
number of points ([GP]).  E. Kussell has recovered the Gromov-Witten
invariants of $\proj^2$ blown up at 2 points by Caporaso and Harris'
``rational fibration method'' ([Ku]).  In another direction, extending work
of I. Vainsencher ([Va]), Kleiman and Piene have examined systems with an
arbitrary, but fixed, number $\de$ of nodes ([K2]).  The postulated number
of $\de$-nodal curves is given (conjecturally) by a polynomial, and they
determine the leading coefficients, which are polynomials in $\de$.
L. G\"{o}ttsche has recently conjectured a surprisingly simple generating
function ([G]) for these polynomials which reproduce the results of
Vainsencher as well as Kleiman and Piene and experimentally reproduce the
numbers of [CH3], [V3], S.T. Yau and E. Zaslow's count of rational curves
on K3-surfaces ([YZ]), and others.  The numbers of curves are expressed in
terms of four universal power series, two of which G\"{o}ttsche gives
explicitly as quasimodular forms.

The philosophy here is that of Caporaso and Harris in [CH3]: we degenerate
the point conditions to lie on $E$ one at a time.  Our perspective,
however, is different: we use the moduli space of stable maps rather
than the Hilbert scheme.

\subsection{Examples} Suppose $d$ and $g$ are integers ($d >0$, $g
\geq 0$), and $s$ is a sequence of non-negative integers with no more
than five terms.  Let $N^g_{d,s}$ be the number of irreducible genus
$g$ degree $d$ plane curves with fixed multiple points with
multiplicities given by $s$, passing through the appropriate number of
general points.  As an example of the algorithm, Table
\ref{table1} gives values of $N^g_{d,s}$ for $d \leq 5$.
When $g=0$, these numbers agree with those found by Pandharipande and
G\"{o}ttsche ([GP] Subsection 5.2; the numbers were called $N_{d,s}$
there).  When $s=0$, these numbers agree with those found by Caporaso
and Harris ([CH3]).

A conjectural algorithm for computing the genus $g$ Gromov-Witten
invariants of the cubic surface (i.e. the plane blown up at six points, no
three collinear and not all on a conic) is given in Subsection
\ref{cubic}.  Based on that conjecture, we compute $N^0_{6,(2^6)}=3240$. 

\begin{table}
\begin{center}
\begin{tabular}{|c|c|c|c|c|c|c|c|c|c|c|c|}
\hline
$N^0_1$ & $N^0_2$ & $N^1_3$ & $N^0_3$ & $N^0_{3,(2)}$ & $N^3_4$ &
$N^2_4$ & $N^1_4$ & $N^0_4$ & $N^2_{4,(2)}$ & $N^1_{4,(2)}$ &
$N^0_{4,(2)}$  \\
\hline
1 & 1 & 1 & 12 & 1 & 1 & 27 & 225 & 620 & 1 & 20 & 96 \\
\hline
\end{tabular}
\end{center}
\begin{center}
\begin{tabular}{|c|c|c|c|c|c|c|c|c|c|c|}
\hline
$N^1_{4,(2^2)}$ & $N^0_{4,(2^2)}$& $N^0_{4,(2^3)}$ & $N^0_{4,(3)}$ & $N^6_5$ & $N^5_5$ & $N^4_5$ &
$N^3_5$ & $N^2_5$ & $N^1_5$ & $N^0_5$ \\
\hline
 1 & 12  & 1 & 1 & 1 & 48 & 882 & 7915 & 36855 & 87192 & 87304 \\
\hline
\end{tabular}
\end{center}
\begin{center}
\begin{tabular}{|c|c|c|c|c|c|}
\hline
$N^5_{5,(2)}$ & $N^4_{5,(2)}$
& $N^3_{5,(2)}$ & $N^2_{5,(2)}$ & $N^1_{5,(2)}$ & $N^0_{5,(2)}$  
 \\
\hline
 1 & 41 & 615 & 4235 & 13775 & 18132   \\
\hline
\end{tabular}
\end{center}
\begin{center}
\begin{tabular}{|c|c|c|c|c|c|c|c|c|}
\hline
$N^4_{5,(2^2)}$ & $N^3_{5,(2^2)}$ & $N^2_{5,(2^2)}$ & $N^1_{5,(2^2)}$
& $N^0_{5,(2^2)}$  & $N^3_{5,(2^3)}$ & $N^2_{5,(2^3)}$ &
$N^1_{5,(2^3)}$ & $N^0_{5,(2^3)}$   \\
\hline
1 & 34 & 396 & 1887 & 3510 &  1 & 27 & 225 & 620   \\
\hline
\end{tabular}
\end{center}
\begin{center}
\begin{tabular}{|c|c|c|c|c|c|c|c|c|}
\hline
$N^2_{5,(2^4)}$ & $N^1_{5,(2^4)}$
& $N^0_{5,(2^4)}$ & $N^1_{5,(2^5)}$ & $N^0_{5,(2^5)}$ & 
$N^3_{5,(3)}$ & $N^2_{5,(3)}$ & $N^1_{5,(3)}$ & $N^0_{5,(3)}$ \\
\hline
1 & 20 & 96 & 1 & 12 & 1 & 28 & 240 & 640 \\
\hline
\end{tabular}
\end{center}
\begin{center}
\begin{tabular}{|c|c|c|c|c|c|c|}
\hline
$N^2_{5,(3,2)}$ & $N^1_{5,(3,2)}$ & $N^0_{5,(3,2)}$ &
$N^1_{5,(3,2^2)}$ & $N^0_{5,(3,2^2)}$ & $N^0_{5,(3,2^3)}$ &
$N^0_{5,(4)}$ \\
\hline
 1 & 20 & 96 & 1 & 12 & 1 & 1 \\
\hline
\end{tabular}
\end{center}
\caption{Numbers of plane curves with fixed multiple points}
\label{table1}
\end{table}

\subsection{Acknowledgements}
The author is grateful to J. Harris, S. Kleiman, R. Pandharipande, and
T. Graber for useful advice and conversations.  This work was developed and
largely written while the author was enjoying the hospitality of the
Mittag-Leffler Institute in May 1997, and he is grateful to the organizers
of the special year in Quantum Cohomology for this opportunity.  This
research was supported by a Sloan Dissertation Fellowship.

\section{Statement of results}

For any sequence $\al = (\al_1, \al_2, \dots)$ of nonnegative integers
with all but finitely many $\al_i$ zero, set 
$$
| \al | = \al_1 + \al_2 + \al_3 + \dots 
$$
$$
I  \al
 = \al_1 + 2\al_2 + 3\al_3 + \dots 
$$
$$
I^\al
 = 1^{\al_1}  2^{\al_2}  3^{\al_3}  \dots 
$$
and
$$
\al ! = \al_1 ! \al_2 ! \al_3! \dots .
$$ We denote by $\lcm(\al)$ the least common multiple of the set $\#
\{ i : \al_i \neq 0 \}$.  The zero sequence will be denoted 0.

We denote by $e_k$ the sequence $(0, \dots, 0, 1, 0 , \dots)$ that is
zero except for a 1 in the $k^{\tth}$ term (so that any sequence $\al
= (\al_1, \al_2, \dots)$ is expressible as $\al = \sum \al_k e_k$).
By the inequality $\al \geq \al'$ we mean $\al_k \geq \al_k'$ for all
$k$; for such a pair of sequences we set $$
\binom \al {\al'} = {\frac { \al!}  {\al' ! (\al - \al')!}} = 
\binom {\al_1} {\al'_1}
\binom {\al_2} {\al'_2}
\binom {\al_3} {\al'_3} \dots .
$$
This notation follows [CH3] and [V3].

Let $H$ be the divisor class of a line in $\proj^2$.  Fix a degree $d$, a
genus $g$, sequences $\al$ and $\be$, and a collection of points $\Ga = \{
p_{i,j} \}_{1 \leq j \leq \al_i}$ (not necessarily distinct) of $E$.  We
define the {\em generalized Severi variety} $\wdg(\al,\be,\Ga)$ to be the
closure (in $|dH|$) of the locus of irreducible reduced curves in $\proj^2$
in class $dH$ of geometric genus $g$, not containing the conic $E$, with
(informally) $\al_k$ ``assigned'' points of contact of order $k$ and
$\be_k$ ``unassigned'' points of contact of order $k$ with $E$.  Formally,
we require that, if $\nu: C^\nu \rightarrow C$ is the normalization of $C$,
then there exist $|\al|$ points $q_{i,j} \in C^\nu$, $j=1, \dots, \al_i$
and $|\be|$ points $r_{i,j} \in C^\nu$, $j=1, \dots, \be_i$ such that 
$$
\nu(q_{i,j}) = p_{i,j} \quad \text{and} \quad \nu^*(E) = \sum i \cdot
q_{i,j} + \sum i \cdot r_{i,j}.
$$
If $I
\al + I\be \neq dH \cdot E = 2d$, $\wdg(\al,\be,\Ga)$ is empty.

For convenience, let 
\begin{eqnarray*}
\Up = \Up^{d,g}(\be) &:=& - (K_{\proj^2} + E) \cdot (dH) +
|\be| + g-1 \\
&=& d + |\be| + g-1.
\end{eqnarray*}
Then $\wdg(\al,\be,\Ga)$ is a projective
variety of pure dimension $\Up$ (Proposition \ref{bigdim}).

{\em Notational warning:} The notation $\vdg(\al,\be,\Ga)$ and
$\wdg(\al,\be,\Ga)$ is used in [CH3] to refer to a slightly different
notion: the sequences $\al$ and $\be$ refer there to tangencies with a
fixed line, rather than the fixed conic $E$.

\begin{defn}
\label{simple} 
If $\{ p_{i,j} \}$ are distinct points except $s_k$ points $\{ p_{1,j}
\}$ are the same $(1 \leq k \leq l)$, we will say that $\Ga$ is {\em
simple}.
\end{defn}
In this case, if $\al'= \al + | s | e_1$,
$\wdg(\al',\be,\Ga)$ generically parametrizes curves that have 
multiple points of order $s_1$, \dots, $s_l$ on $E$ (with each branch
transverse to $E$), $i$-fold tangent to $E$ at $\al_i$ fixed points of
$E$, and $i$-fold tangent to $E$ at $\be_i$ other points of $E$.  When
discussing properties of $\wdg(\al',\be,\Ga)$ that depend only on
$\al$, $\be$, and $s$, we write $W^{d,g}(\al,\be,s)$ for convenience.

Let $N_{\ti}^{d,g}\abG$ be the number of points of $\wdg\abG$ whose
corresponding curve passes through $\Up$ fixed general points of $\proj^2$.
If $\Ga$ is {\em simple} then $N^{d,g}_{\ti}(\al',\be,\Ga)$ depends only on
$(\al,\be,s)$ (Subsection \ref{therecursiveformulas}), so we write
$N^{d,g}_{\ti}\abs$.  Then $N_{\ti}^{d,g}\abs$ is the degree of the
generalized Severi variety (in the projective space $|dH|$).  The main
result of this note is the following.

\begin{tm}
\label{irecursion}
If $\dim \wdg\abs>0$, then 
\begin{eqnarray*}
N_{\ti}^{d,g}\abs &=&  \sum_{\beta_k > 0} k N_{\ti}^{d,g}(\alpha + e_k,
\beta - e_k,s)
\\
& & + \sum \frac 1 \si \binom {\al} {\al^1, \dots, \al^l, \al- \sum \al^i} \binom {\Up^{d,g}(\be)-1} {\Up^{d^1,g^1}(\be^1), \dots, \Up^{d^l,g^l}(\be^l)} \\
& & \cdot  \prod_{i=1}^l \binom {\be^i} {\ga^i} I^{\be^i - \ga^i} N_{\ti}^{d^i,g^i}(\al^i,\be^i,s^i) 
\end{eqnarray*}
where the second sum runs over choices of $d^i, g^i, \alpha^i,
\beta^i, \gamma^i, s^i$ ($1 \le i \le l$), where $d^i$ is a positive integer,
$g^i$ is a non-negative integer, $\alpha^i$, $\be^i$, $\ga^i$, $s^i$ are
sequences of non-negative integers, $\sum_i d^i = d-2$, $\sum_i \ga^i
= \be$, $\be^i \gneq \ga^i$, $\sum_i s^i_k = s_k - 1$ or $s_k$, and $\si$ is the number of
symmetries of the set $\{ (d^i,g^i,\al^i,\be^i,\ga^i, s^i) \}_{1 \leq i
\leq l}$.
\end{tm}
In the second sum, for the summand to be non-zero, one must also have
$\sum \al^i \leq \al$, and $I \al^i + I \be^i + |s^i|=
2 d^i$.  

If $\tilde{s}$ is the sequence $s$ with the zeros removed, then clearly 
\begin{equation}
\label{fred}
\wdg\abs = \wdg(\al,\be,\tilde{s})
\end{equation}
If $a$ is the number of ones in $s$, and $\tilde{s}$ is the sequence $s$ with the ones removed, then clearly
\begin{equation}
\label{barney}
\wdg\abs = \wdg(\al + a e_1,\be,\tilde{s}).
\end{equation}
(Requiring a curve to have a multiplicity-1 multiple point at a fixed point
of $E$ is the same as requiring the curve to pass through a fixed point of
$E$.)  If $s_k>1$ for all $k$, then the variety $\wdg\abs$ has dimension 0
if and only if $d=1$, $g=0$, $\be=0$, $s=0$, and $\al = 2 e_1$ or $e_2$ (from
Proposition \ref{bigdim} and simple case-checking).  The first case is a line
through two fixed points of the conic $E$, and the second is a line tangent
to $E$ at a fixed point.  In both cases, $N^{d,g}_{\ti}\abs = 1$.  Therefore,
with this ``seed data'', Theorem \ref{irecursion} provides a means of
recursively computing $N_{\ti}^{d,g}\abs$ for all $d$, $g$, $\al$, $\be$,
$s$.

\subsection{Relationship to Gromov-Witten invariants}

Let $B$ be the blow-up of the plane at $l$ points ($0 \leq l \leq 5$),
no three collinear.  Let $E$ be a smooth conic through the $l$ points
(unique if $l=5$).  Let $H$ be the pullback of a line to $B$, and
$E_1$, \dots, $E_l$ the exceptional divisors.  The irreducible genus
$g$ curves in a divisor class $D$ on $B$, through an appropriate
number ($-K_B \cdot D + g-1$) of points, can now be counted.  Call
this number $GW^{D,g}_B$ for convenience.  If $D$ is an exceptional
divisor, or the proper transform of $E$ if $l=5$, then 
$$
GW^{D,g}_B =
\begin{cases}
1 & \text{if $g=0$,} \\
0 & \text{otherwise.}
\end{cases}
$$
If $D= dH - \sum_{k=1}^l m_k E_k$ is any
other divisor class, then the (finite number of) genus $g$ curves in
class $D$ on $B$ through $-K_B \cdot D + g-1$ general points
correspond to the degree $d$ plane curves with $l$ fixed multiple points
of multiplicity $m_1$, \dots, $m_5$ through the same number of points.
Thus if $s_i = D \cdot E_i$, then
$$
GW^{D,g}_B = N^{d,g}(0,(2d-\sum m_i) e_1, s).
$$

In Subsection \ref{gwenumerative}, it will be shown that the numbers
$GW^{D,g}_B$ give all the genus $g$ Gromov-Witten invariants of $B$.  (It
was previously known that the genus 0 invariants are enumerative, see [RT]
Section 10 and [GP] Lemma 4.10.)

\subsection{The strategy}

In order to understand generalized Severi varieties, we will analyze
certain moduli spaces of maps.  Let $\mbar_g(\proj^2,d)$ be the moduli
space of maps $\pi: C \rightarrow \proj^2$ where $C$ is irreducible,
complete, reduced, and nodal, $(C, \pi)$ has finite automorphism
group, and $\pi_* [C] = dH$.  Let $d$, $g$, $\al$, $\be$, $\Ga$ be as
in the definition of $\wdg(\al,\be,\Ga)$ above.  Define the {\em
generalized Severi variety of maps} $\wdg_m(\al,\be,\Ga)$ to be the
closure in $\mbar_g(\proj^2,d)$ of points representing maps $(C,\pi)$ where
each component of $C$ maps birationally to its image in $\proj^2$, no
component maps to $E$, and $C$ has (informally) $\al_k$ ``assigned''
points of contact of order $k$ and $\be_k$ ``unassigned'' points of
contact of order $k$ with $E$.  Formally, we require that there exist
$|\al|$ smooth points $q_{i,j} \in C$, $j=1, \dots, \al_i$ and $|\be|$
smooth points $r_{i,j} \in C$, $j=1, \dots, \be_i$ such that
$$
\pi(q_{i,j}) = p_{i,j} \quad \text{and} \quad \pi^*(E) = \sum i \cdot q_{i,j} +
\sum i \cdot r_{i,j}.
$$

There is a natural rational map from each component of
$\wdg\abs$ to $\wdg_m\abs$, and the 
dimension of the image will be $\Up$.  We will prove:

\begin{pr}
\label{idim}
The components of $\wdg_m\abs$ have dimension at most $\Up$,
and the union of those with dimension exactly $\Up$ is
the closure of the image of $\wdg\abs$ in $\wdg_m\abG$.
\end{pr}

(This will be an immediate consequence of Theorem \ref{bigdim}.)

Fix $\Up$ general points $s_1$, \dots, $s_\Up$ on $\proj^2$.  The image of
the maps in $\wdg_m\abs$ whose images pass through these points are
reduced.  ({\em Proof:} Without loss of generality, restrict to the
union $W$ of those components of $\wdg_m\abs$ with dimension $\Up$.  By
Proposition \ref{idim}, the subvariety of $W$ corresponding to maps
whose images are {\em not} reduced contains no components of
$W$ and hence has dimension less than $\Up$.  Thus no image of such a
map passes through $s_1$, \dots, $s_{\Up}$.)

Therefore, if $H$ is the divisor class on $\wdg_m\abs$ corresponding to
requiring the image curve to pass through a fixed point of $\proj^2$, then
$$
N_{\ti}^{d,g}\ab = H^\Up.  
$$ 

Define the {\em intersection dimension} of a family $V$ of maps to
$\proj^2$ (denoted $\idim V$) as the maximum number $n$ of general
points $s_1$, \dots $s_n$ on $\proj^2$ such that there is a map $\pi: C
\rightarrow \proj^2$ in $V$ with $\{ s_1, \dots, s_n \} \subset
\pi(C)$.  Clearly $\idim V \leq \dim V$.

Our strategy is as follows.  Fix a general point $q$ of $E$.  Let
$H_q$ be the Weil divisor on $\wdg_m\abs$ corresponding to maps with images
containing $q$.  We will find the components of $\wdg_m\abs$ with
intersection dimension $\Up-1$ and relate them to
$W^{d',g'}_m(\al',\be',s')$ for appropriately chosen $d'$, $g'$, $\al'$,
$\be'$, $s'$.  Then we compute the multiplicity with which each of these
components appears.  Finally, we derive a
recursive formula for $N_{\ti}^{d,g}\abs$ (Theorem \ref{irecursion}).

\subsection{Counting reducible curves}

Analogous definitions can be made of the space $\vdg_m\abs$
parametrizing possibly reducible curves.  (Then $\wdg_m\abs$ is the
union of connected components of $\vdg_m\abs$ where the source curve
$C$ is connected.)  The arguments in this case are identical,
resulting in a recursive formula for $N^{d,g}\abs$, the number of maps
$(C,\pi)$ from genus $g$ curves with $\pi_*[C]=dH$, and intersection
of $\pi_* C$ with $E$ determined by $\al$ and $\be$, passing through
$\Up$ fixed general points of $\proj^2$:

\begin{tm}
\label{rrecursion}
If $\dim \vdg_m\abs>0$, then 
\begin{eqnarray*}
N^{d,g}\abs = \sum_{\be_k > 0} k N^{d,g}(\al + e_k, \be-e_k,s) 
\\
+ \sum I^{\be'-\be} {\binom \al {\al'}} \binom {\be'}{\be} 
N^{d-2,g'}(\al',\be',s')
\end{eqnarray*}
where the second sum is taken over all $\al'$, $\be'$, $g'$ satisfying
$\al' \leq \al$, $\be' \geq \be$, $s'_k= s_k$ or $s_k-1$, $g-g' =
|\be'-\be| - 1$, $I \al' + I \be' +|s'| = 2d-4$.
\end{tm}
Although the recurrence is simpler than in Theorem \ref{irecursion}, the
seed data is more complicated.  If $\dim \vdg_m\abs=0$ and $N^{d,g}\abs
\neq 0$, then $\be = 0$, $\al_k = 0$ for $k>2$, $g=1-2d$, and $2d = |s| + I
\al$.  In this case, each point of $\vdg_m\abs$ corresponds to a union of
lines in the plane.  Then if $s'$ is the sequence $s$ with $\al_1$ ones
appended, $N^{d,g}\abs$ is the number of ways of expressing $s'$ as a sum
of sequences that are 0 except for 1's in two places.  (This involves
chasing through the definitions: we are counting maps of $2d$ $\proj^1$'s
to the plane, each of degree 1, intersecting $E$ in a certain way.)
Alternatively, $N^{d,g}\abs$ is the number of graphs with $|s'|$ labeled
vertices with vertex $i$ having valence $s'_i$, where each edge joints two
different vertices (although more than one edge can connect the same pair
of vertices).

\subsection{Simplifying formulas}

Computationally, it is simpler to deal with the formulas of Theorems
\ref{irecursion} and \ref{rrecursion} when $s=0$.  (This was the form
proved in [V3].)  The following lemma allows us to quickly reduce to the
case $s=0$ when $s_i=0$ for $i>3$.

\begin{lm}
\begin{enumerate}
\item[(0)] $N^{d,g}_{\ti}(\al,\be,s \cup (0)) = N^{d,g}_{\ti}\abs$
\item[(1)] $N^{d,g}_{\ti}(\al,\be,s\cup (1)) = N^{d,g}_{\ti}(\al+e_1,\be,s)$
\item[(2)] $N^{d,g}_{\ti}(\al,\be,s \cup (2)) = \frac 1 2 \left(
N^{d,g}_{\ti}(\al + 2 e_1,\be,s) - N^{d,g}_{\ti}(\al + e_2,\be,s) \right)$
\item[(3)] $N^{d,g}_{\ti}(\al,\be,s\cup (3)) = \frac 1 6 
N^{d,g}_{\ti}(\al + 3 e_1,\be,s) - \frac 1 2 N^{d,g}_{\ti}(\al + e_1 +
e_2,\be,s) + \frac 1 3 N^{d,g}_{\ti}(\al + e_3,\be,s)$
\end{enumerate}
\end{lm}
Parts (0) and (1) are tautological, and were stated earlier (equations
(\ref{fred}) and (\ref{barney})).  There are analogous expressions for
$N^{d,g}_{\ti}(\al,\be,s \cup (n))$ for all $n$.  The result still holds
when $N^{d,g}_{\ti}$ is replaced by $N^{d,g}$.  The lemma can be proved by
induction (on $\al$, $\be$, $s$) and Theorem \ref{irecursion}.  It can also
be proven by degeneration methods, and by the study of divisors on
generalized Severi varieties.

\subsection{Variations on a theme}
The same arguments provide means of computing the number of curves in
$V^{d,g}\abG$ through an appropriate number of points for $(X,E) =
(\proj^2,H)$, ($\proj^2,E)$, or $(\eff_n,E)$, even if $\Ga$ is not a
collection of reduced points.

\section{Proof of results}

\subsection{Dimension counts}

The pair $(X,E) = (\proj^2,E)$ satisfies properties P1--P4 of [V3],
Subsection 1.1:
\begin{enumerate}
\item[{\bf P1.}]  $X$ is a smooth surface and $E \cong \proj^1$ is a divisor on $X$.
\item[{\bf P2.}]
The surface $X \setminus E$ is minimal, i.e. contains no (-1)-curves.
\item[{\bf P3.}]
The divisor class $K_X + E$ is negative on every curve on $X$.
\item[{\bf P4.}] 
If $D$ is an effective divisor such that $-(K_X + E) \cdot D = 1$,
then $D$ is smooth.
\end{enumerate}

Recall [V3] Theorem 2.1, whose proof depended only on
properties P1--P4:

\begin{tm}
\label{bigdim}
\begin{enumerate}
\item[(a)]  Each component $V$ of $\vdg_m\abG$ is of dimension at most 
$$
\Up = \Up^{d,g}(\be) = -(K_X + E) \cdot (dH) + | \be | + g-1.
$$
\item[(b)]
The stable map $(C, \pi)$ corresponding to a general point of any
component of dimension $\Up$ satisfies the following properties.
\begin{enumerate}
\item[(i)] The curve $C$ is smooth, and the map $\pi$ is an immersion.
\item[(ii)] The image is a reduced curve.  If $\Ga$ consists of distinct points, then the image is smooth along its intersection with $E$.
\end{enumerate}
\item[(c)]  Conversely, any component whose general map satisfies 
property (i) has dimension $\Up$.
\end{enumerate}
\end{tm}

By ``the image is a reduced curve'', we mean $\pi_*[C]$ is a sum of
distinct irreducible divisors on $X$.

\subsection{Determining intersection components and their multiplicities}

Fix $d$, $g$, $\al$, $\be$, $s$, $\Ga$, and a general point $q$ on
$E$.  Let $H_q$ be the divisor on $\vdg_m\abs$ corresponding to maps
whose image contains $q$.  We will derive a list of subvarieties (which
we will call {\em intersection components}) in which each component of
$H_q$ of intersection dimension $\Up-1$ appears, and then calculate
the multiplicity of $H_q$ along each such component.

In [V3], in a more general situation, a list of intersection
components was derived and the multiplicities calculated.  We recall
these results, and apply them in this particular case.

The potential components come in two classes that
naturally arise from requiring the curve to pass through $q$.  First, one of
the ``moving tangencies'' $r_{i,j}$ could map to $q$.  We will call
such components {\it Type I intersection components}.

Second, the curve could degenerate to contain $E$ as a component.  We
will call such components {\it Type II intersection components}.  For
any sequences $\al'' \leq \al$, $\ga \geq 0$, and subsets $\{
p''_{i,1}, \dots, p''_{i,\al''_i} \}$ of $\{ p_{i,1}, \dots,
p_{i,\al_i} \}$, let $g'' = g + |\ga| + 1$ and $\Ga'' = \{ p''_{i,j}
\}_{1 \leq j \leq \al''_i}$.  Define 
$K(\al'',\be,\ga,\Ga'')$ as the closure in $\mbar_g(X,d)'$ of points
representing maps $\pi: C' \cup C''
\rightarrow X$ where
\begin{enumerate}
\item[K1.] the curve $C'$ maps isomorphically to $E$,
\item[K2.] the curve $C''$ is smooth, each component of $C''$ maps
birationally  
to its image, no 
component of $C''$ maps to $E$, and there exist $|\al''|$ points
$q_{i,j} \in C''$, $j = 1$, \dots, $\al_i''$, $|\be|$ points $r_{i,j}
\in C''$, $j = 1$, \dots, $\be_i$, $|\ga|$ points $t_{i,j} \in C''$,
$j = 1$, \dots, $\ga_i$ such that $$
\pi(q_{i,j}) = p''_{i,j} \quad \text{and} \quad
(\pi|_{C''})^*(E) = \sum i \cdot q_{i,j} + 
 \sum i \cdot r_{i,j} + 
 \sum i \cdot t_{i,j},
$$
and
\item[K3.] the intersection of the curves $C'$ and $C''$ is $\{ t_{i,j} \}_{i,j}$.
\end{enumerate}

The variety $K(\al'',\be,\ga,\Ga'')$ is empty unless $I(\al''+\be+\ga)
= (dH-E) \cdot E = 2(d-2)$.  The genus of $C''$ is $g''$, and there is a degree
$\binom {\be+\ga} \be$ rational map
\begin{equation}
\label{rratmap}
K(\al'',\be,\ga,\Ga'') \dashrightarrow V_m^{dH-E,g''}(\al'',\be+\ga,\Ga'')
\end{equation}
corresponding to ``forgetting the curve $C'$''.

Define the Type II intersection component $K_V(\al'',\be,\ga,\Ga'')$
in the same way as $K(\al'',\be,\ga,\Ga'')$, with an additional
condition:
\begin{enumerate}
\item[K4.] The collection $\Ga' = \Ga \setminus \Ga''$ consists of
distinct points.
\end{enumerate}
If for the general $(C,\pi)$ in $\vdg_m\abs$, $\pi(C)$ has a $k$-fold
point at some fixed $p$ in $E$, then for a general $(C' \cup C'',
\pi)$ in $K_V(\al'',\be,\ga,\Ga'')$, $\pi(C'')$ has at least a $(k-1)$-fold point at $p$,
by condition K4.

\begin{tm}
\label{rlist}
Fix $d$, $g$, $\al$, $\be$, $\Ga$, and a point $q$ on $E$ not in $\Ga$.
Let $K$ be an irreducible component of $H_q$ with 
intersection dimension $\Up  - 1$.  Then set-theoretically, either
\begin{enumerate}
\item[I.] $K$ is a component of $\vdg_m(\al+e_k,\be-e_k,\Ga')$, 
where $\Ga'$ is the same as $\Ga$ except $p'_{k,\al_{k+1}} = q$, or
\item[II.]  $K$ is a component of $K_V(\al'',\be,\ga,\Ga'')$ for some $\al''$, $\ga$, $\Ga''$.
\end{enumerate}
\end{tm}

\bpf
By [V3], Theorem 3.1, either
\begin{enumerate}
\item[I.] $K$ is a component of $\vdg_m(\al+e_k,\be-e_k,\Ga')$, 
where $\Ga'$ is the same as $\Ga$ except $p'_{k,\al_{k+1}} = q$, or
\item[II.]  $K$ is a component of $K(\al'',\be,\ga,\Ga'')$ for some $\al''$, $\ga$, $\Ga''$.
\end{enumerate}

Suppose the point $p$ appears in $\Ga$ $n$ times.  Then for a general
map $(C,\pi)$ in $\vdg_m\abG$, $\pi^*(p)$ is a length $n$ scheme by
Theorem \ref{bigdim}(b) above.  If $K$ is a component of
$K(\al'',\be,\ga,\Ga'')$ for some $\al''$, $\ga$, and $\Ga''$, and
$(C,\pi) = (C' \cup C'', \pi)$ is a general map in $K$, then $(\pi
|_{C'})^* p$ has length 1, so $(\pi |_{C''})^* p$ must have length at
least $n-1$.  As $\pi|_{C''}$ is an immersion (by Theorem
\ref{bigdim}(b)(i)), $(\pi|_{C''})^* p$ is the number of times $p$
appears in $\Ga''$.  Thus if $p$ appears $n$ times in $\Ga$ then it
appears at least $n-1$ times in $\Ga''$, so it appears at most once in
$\Ga' = \Ga \setminus \Ga''$.  Hence $K$ is actually a component of
$K_V(\al'',\be,\ga,\Ga'')$.  \epf

There are other components of the divisor $H_q$ not counted in
Theorem \ref{rlist}, but they will be enumeratively
irrelevant.  (See the end of [V3] Section 3 for example of such
behavior.)

There is an analogous result for $\wdg_m\abs$.  Let 
$K_W(\al'',\be,\ga,\Ga'')$ be the union of components of 
$K_V(\al'',\be,\ga,\Ga'')$ where the source curve $C$ is connected.  Then $C''$ is a union of $l$ curves with image of degree $d^i$, or arithmetic genus $g^i$, with a subset $\Ga^i$ of the points $\Ga$ (and induced sequence $\al^i \leq \al$), and induced 
sequences $\be^i$, $\ga^i$, $s^i$.  As the source curve is connected, $\ga^i>0$ for all $i$.

Let $\si$ be the symmetry group of the data
$(d^i,g^i,\al^i,\be^i,\ga^i,s^i)$ (so for example $\si$ is the
one-element group if no two $(d^i,g^i,\al^i,\be^i,\ga^i,s^i)$ are the
same, and $\si = S_l$ is they are all the same).  Then there is a
generically $|\si|$-to-1 cover of this component $K_W$ of
$K_W(\al'',\be,\ga,\Ga'')$ which distinguishes $C^1$, \dots, $C^l$.
Then there is a degree $\prod \binom { \be^i + \ga^i} {\ga^i}$ map
\begin{equation}
\label{iratmap}
K_V \dashrightarrow \prod V^{d^i,g^i} (\al^i, \be^i + \ga^i, \Ga^i).
\end{equation}

\begin{tm}
\label{ilist}
Fix $d$, $g$, $\al$, $\be$, $\Ga$, and a point $q$ on $E$ not in $\Ga$.
Let $K$ be an irreducible component of $H_q$ on $\wdg_m\abG$ with 
intersection dimension $\Up  - 1$.  Then set-theoretically, either
\begin{enumerate}
\item[I.] $K$ is a component of $\wdg_m(\al+e_k,\be-e_k,\Ga')$, 
where $\Ga'$ is the same as $\Ga$ except $p'_{k,\al_{k+1}} = q$, or
\item[II.]  $K$ is a component of $K_W(\al'',\be,\ga,\Ga'')$ for some $\al''$, $\ga$, $\Ga''$.
\end{enumerate}
\end{tm}

\bpf 
As $\wdg_m\abG$ is the union of components of $\vdg_m\abG$ where
the source is connected, this follows immediately from Theorem
\ref{rlist}.  \epf

\subsection{Multiplicity of $H_q$ along intersection components}

The proofs of the multiplicity calculations are identical to those of [V3].

Let $K_k$ be the union of Type I intersection components of the form
$\vdg_m(\al+e_k, \be-e_k, \Ga')$ as described in Theorem \ref{rlist}.
The following proposition is [V3] Proposition 4.1.

\begin{pr}
\label{multI}
The multiplicity of $H_q$ along $K_k$ is $k$.
\end{pr}

Suppose $K = K_V(\al'', \be, \ga, \Ga'')$ is a Type II component of
$H_q$ (on $\vdg_m\ab$).  Let $m_1$, \dots, $m_{|\ga|}$ be a set of
positive integers with $j$ appearing $\ga_j$ times ($j = 1$, 2,
\dots), so $\sum m_i = I \ga$.  The following proposition
is [V3] Proposition 5.2.

\begin{pr}
\label{multII}
The multiplicity of $H_q$ along $K$ is $m_1 \dots m_{|\ga|} = I^{\ga}$.
\end{pr}

(The proof in [V3] assumed only that $\Ga \setminus \Ga''$ consisted
of distinct points.)

As $\wdg_m\abG$ is a connected union of components of $\vdg_m\abG$,
analogous multiplicity results hold for $\wdg_m\abG$.

\subsection{The Recursive Formulas}
\label{therecursiveformulas}
Let $H_q$ be the
divisor on $\vdg_m(\al,\be,\Ga)$ corresponding to requiring the image
to contain a general point $q$ of $E$.  The components of $H_q$ of
intersection dimension $\Up - 1$ were determined in Theorem
\ref{rlist}, and the multiplicities were determined in Propositions
\ref{multI} and \ref{multII}: 
\begin{pr}  
In the Chow ring of $\vdg_m\abG$, modulo Weil divisors of intersection
dimension less than $\Up - 1$,
$$
H_q = \sum_{\be_k>0} k \cdot \vdg_m(\al+e_k,\be-e_k,\Ga \cup \{ q \}) 
+ \sum I^{\ga} \cdot K_V( \al'', \be, \ga, \Ga'')
$$
where the second sum is over all $\al'' \leq \al$, $\ga \geq 0$, $\Ga'' = \{
p''_{i,j} \}_{1 \leq j \leq \al''_i} \subset \Ga$ such that $\Ga \setminus \Ga''$ consists of distinct points, 
$I(\al'' + \be + \ga ) = (dH) \cdot E$. 
\end{pr}
Intersect both sides of the equation with $H_q^{\Up - 1}$.  As those
dimension $\Up - 1$ classes of intersection dimension less than $\Up -
1$ are annihilated by $H_q^{\Up - 1}$, we still have equality:
\begin{eqnarray*}
N^{d,g}\abG &=& H_q^{\Up} \\
&=&
\sum_{\be_k>0} k \vdg_m(\al+e_k,\be-e_k,\Ga \cup \{ q \}) \cdot H_q^{\Up - 1}\\
& & + \sum I^{\ga} \cdot K_V( \al'', \be, \ga, \Ga'') \cdot H_q^{\Up - 1}.  
\end{eqnarray*}
From (\ref{rratmap}), each $K_V(\al'', \be, \ga, \Ga'')$
admits a degree $\binom {\be + \ga}
\be$ rational map to $V^{dH-E,g''}_m(\al'', \be+\ga,\Ga'')$ 
(where $g'' = g- |\ga| + 1$) corresponding to ``forgetting the component mapping to $E$'', so
$$
K_V(\al'', \be, \ga, \Ga'') \cdot H_q^{\Up - 1} = 
\binom {\be + \ga} \ga N^{d-2,g''}(\al'', \be+\ga).
$$
For each $\al''$, there are $\binom \al {\al''}$ choices of $\Ga''$
(as this is the number of ways of choosing 
$\{ p''_{i,1} , \dots, p''_{i,\al''_i} \}$ from
$\{ p_{i,1} , \dots, p_{i,\al_i} \}$).  Thus
\begin{eqnarray}
\nonumber
N^{d,g}\abG &=& \sum_{\be_k > 0} k N^{d,g}(\al + e_k, \be-e_k, \Ga \cup \{
q \}) 
\\
& & + \sum I^{\ga} {\binom \al {\al''}} \binom {\be + \ga}{\be} 
N^{d-2,g''}(\al'',\be + \ga, \Ga'').
\label{donaldduck}
\end{eqnarray}
If $\Ga$ is {\em simple} (Definition \ref{simple}), then so are $\Ga \cup
\{ q \}$ (as $q$ is general) and $\Ga''$ (as $\Ga'' \subset \Ga$).  By
(\ref{donaldduck}) and induction, if $\Ga$ is simple,
$N^{d,g}(\al',\be,\Ga)$ depends only on $(\al, \be, s)$.  Rewriting
(\ref{donaldduck}) in terms of $\al$, $\be$, and $s$, this is Theorem
\ref{rrecursion}.

By the same argument for the irreducible case, using the rational map (\ref{iratmap}) rather than (\ref{rratmap}), yields Theorem \ref{irecursion}.

\subsection{Genus $g$ Gromov-Witten invariants are enumerative on Fano
surfaces}
\label{gwenumerative}

The definition of (genus $g$) Gromov-Witten invariants $I_{g,D}(\ga_1
\cdots \ga_n)$, where $g$ is the genus, $D \in A_1 X$, and $\ga_i \in A^*
X$ is given in [KM] and summarized in [V3].  As in [V3], to compute
Gromov-Witten invariants of a surface, it suffices to deal with the case
where $D$ is effective and non-zero and the $\ga_i$ are points.  Recall
[V3] Lemma 7.1:
\begin{lm}
\label{gwlemma}
Let $X$ be a Fano surface, and let $D$ be an effective divisor class
on $X$.  Suppose that $M$ is an irreducible component of $\mbar_g(X,d)$
with general map $(C,\pi)$.  Then
$$
\idim M \leq -K_X \cdot D + g-1.
$$
If equality holds and $D \neq 0$, then $\pi$ is an immersion.
\end{lm}

Thus the genus $g$ Gromov-Witten invariants of $\fn$ can be computed as
follows.  We need only compute $I_{g,D}(\ga_1
\cdots \ga_n)$ where $D$ is effective and nonzero, and the $\ga_i$ are
(general) points.  By Lemma \ref{gwlemma}, this is the number of immersed
genus $g$ curves in class $D'$ through the appropriate number of points of
$X$.  If $D$ is a (-1)-curve, the number is 1.  Otherwise, the number is
recursively calculated by Theorem \ref{irecursion}.

\subsection{The cubic surface}
\label{cubic}
By the previous subsection, the algorithm of Theorem \ref{irecursion}
computes the genus $g$ Gromov-Witten invariants of the plane blown up at up
to five points.  As GW-invariants are deformation-invariant, one might hope
to compute the invariants of the plane blown up at six points $P_1$, \dots,
$P_6$ in general position (i.e. a general cubic surface in $\proj^3$) by
degenerating the six points $P_1$, \dots, $P_6$ to lie on a conic $E$.
Call the resulting surface $B'$.  (Then $B'$ has a (-2)-curve, the proper
transform of $E$.  The canonical map is an embedding away from $E$, and $E$
is collapsed to a simple double point.)  If the enumerative significance of
the genus $g$ Gromov-Witten invariants on this surface could be determined,
Theorem \ref{irecursion} could be used to determine the invariants of
$B'$, and hence (by deformation-invariance of Gromov-Witten invariants) the
invariants of any cubic hypersurface.  

The rational ruled surface $\eff_2$ ($\proj( \oh_{\proj^1} \oplus
\oh_{\proj^1}(2))$) is similar, in that $\eff_2$ has a (-2)-curve, and the
canonical map is an embedding away from the (-2)-curve, which
is collapsed to a simple double point.  In [K1], p. 22-23, Kleiman gives an
enumerative interpretation for a particular genus 0 GW-invariant of $\eff_2$, which was
explained to him by Abramovich.  This interpretation suggests the following
conjecture.

\begin{conj}
\label{conj}
Suppose $X$ is $\eff_2$ or $B'$, and $E$ is the (-2)-curve on $X$.  Let $D$
be an effective divisor class on $X$ (not 0 or $E$), and $\ga$ the class of
a point.  Then the Gromov-Witten invariant $I_{g,D}(\ga^{- K_X \cdot D +
g-1})$ is the number of maps $\pi:  C \rightarrow X$ with $\pi_*[C] = D$,
where
\begin{enumerate}
\item[(i)] the curve $C$ has one component $C_0$ {\em not} mapping to $E$,
and
\item[(ii)] any other component $C'$ of $C$ maps isomorphically to $E$, and
$C'$ intersections $\overline{C \setminus C'}$ at one point, which is
contained in $C_0$.
\end{enumerate}
\end{conj}
Simple tests on both $\eff_2$ and $B'$ seem to corroborate this
conjecture.  As a more complicated test case, we compute the number of
rational sextic curves in the plane with six nodes at fixed points $P_1$,
\dots, $P_6$, and passing through five other fixed points $Q_1$, \dots,
$Q_5$, where all the points are in general position.  (This is the
Gromov-Witten invariant $N^0_{6,2^6}$ of the cubic surface.)
[DI] p. 119 gives this number as 2376, while [GP] p. 25 gives the number as
3240.  G\"{o}ttsche and Pandharipande checked their number using different
recursive strategies.

According to the conjecture, this invariant is the sum of three
contributions.
\begin{enumerate}
\item Those (irreducible) rational sextics with six fixed nodes $P_1$, \dots,
$P_6$ lying on a conic, passing through $Q_1$, \dots, $Q_5$.  By Theorem
\ref{irecursion} (and some computation), this number is 2002.
\item A stable map $\pi: C \rightarrow \proj^2$ where $C$ has two
irreducible rational components $C_0$ and $C_1$ joined at one point, $\pi$
maps $C_1$ isomorphically to $E$, and $\pi$ maps $C_0$ to an irreducible
rational quartic through $P_1$, \dots, $P_6$ (which lie on a conic) and
$Q_1$, \dots, $Q_5$.  The image of the node $C_0 \cap C_1$ is one of the
two points $\pi(C_0) \cap E \setminus \{ P_1, \dots, P_6 \}$.  By Theorem
\ref{irecursion}, there are 616 such quartics.  There are two choices for
the image of the node $C_0 \cap C_1$, so the contribution is 1232.
\item A stable map $\pi:  C \rightarrow \proj^2$ where $C$ has three
irreducible rational components $C_0$, $C_1$, $C_2$, where $C_1$ and $C_2$
intersect $C_0$, $\pi$ maps $C_1$ and $C_2$ isomorphically to $E$, and
$\pi$ maps $C_0$ isomorphically
to the conic through $Q_1$, \dots, $Q_5$.  There are 12 choices of pairs of
images of the nodes $C_0 \cap C_1$ and $C_0 \cap C_2$, and we must divide
by 2 as exchanging $C_1$ and $C_2$ preserves the stable map.  This this
contribution is 6.
\end{enumerate}
Therefore, assuming Conjecture \ref{conj}, 
$$
N^0_{6,2^6} = 2002 + 1232 + 6 = 3240,
$$
in agreement with [GP].

\noindent
\address{Department of Mathematics \\ Princeton University \\
Fine Hall, Washington Road \\ Princeton, NJ 08544-1000 \\
vakil@math.princeton.edu} 
\end{document}